# Price systems for markets with transaction costs and control problems for some finance problems

Tzuu-Shuh Chiang[1], Shang-Yuan Shiu[2] and Shuenn-Jyi Sheu[1],*

*Academia Sinica and University of Utah*

**Abstract:** In a market with transaction costs, the price of a derivative can be expressed in terms of (preconsistent) price systems (after Kusuoka (1995)). In this paper, we consider a market with binomial model for stock price and discuss how to generate the price systems. From this, the price formula of a derivative can be reformulated as a stochastic control problem. Then the dynamic programming approach can be used to calculate the price. We also discuss optimization of expected utility using price systems.

## 1. Introduction

Duality approach is frequently used for financial problems in incomplete markets. This approach can also be applied to markets with transaction costs. In [12], a discrete market with transaction costs is considered. In the market studied there is a stock and a bond that we can trade. Let $\lambda_1, \lambda_0 > 0$ be the proportional costs for selling and buying the stock. Then the replication cost at time 0 for a portfolio $Y = (Y^0, Y^1)$ at time $T$ is given by

$$\pi^*(Y) = \sup\{E[Y^0 \rho^0 + Y^1 \rho^1]\}. \tag{1.1}$$

The supremum is taken over $(\rho^0, \rho^1)$ ((preconsistent) price systems) which depend on $\lambda_0, \lambda_1$. This will be described in details below.

A similar result for diffusion models is given in [3].

Our interest is to use price systems to calculate the price of a derivative and find optimal strategy for hedging problem. We will also discuss the use of price systems to study portfolio optimization problem. There is a similarity between these problems that they can be reformualted as optimization problems. We shall consider binomial model (it can also be extended to multinomial model) and find a dynamics to generate the price systems $(\rho^0, \rho^1)$. A price system becomes a controlled process. The optimization problems become stochastic control problems. Then dynamic programming approach can be used.

The paper is organized as follows. In Section 2, we give notations and give the framework. In Section 3, we describe price systems and give a price formula for derivatives in terms of price systems. In Section 4, we discuss the optimization of expected utility using price systems. In Sections 5, 6 and 7, we consider binomial models. We present a dynamics to generate the price systems. We reformulate some finance problems as stochastic control problems. Then we use dynamic programming to calculate the value functions.

*Supported by the grant NSC 93-2115-M-001-010.
[1]Institute of Mathematics, Academia Sinica, Nankang, Taipei, Taiwan, e-mail: matsch@math.sinica.edu.tw; sheusj@math.sinica.edu.tw
[2]Department of Mathematics, University of Utah, 155 South 1400 East, Salt Lake City, UT 84112-0090, USA.

## 2. Finite market with one stock

The framework can be described as follows.
We consider $(\Omega, \mathcal{F}, P)$ a finite probability space and $\{\mathcal{F}_k\}$ a filtration. Let

$$P^0(k;\omega), P^1(k;\omega)$$

be the prices for bond and stock. Then $P^0, P^1$ are adapted to $\{\mathcal{F}_k\}$. Define

$$\hat{P}(k;\omega) = P^1(k,\omega)/P^0(k;\omega),$$

the discounted price.

A trading strategy is given by $\{I(k;\omega)\}_{k=0}^T$, a stochastic process adapted to $\{\mathcal{F}_k\}$. $I(k;\omega)$ is the number of shares that the stock is bought or sold,

$$\begin{aligned} I(k;\omega) \geq 0, & \quad \text{buy stock at } k, \\ I(k;\omega) < 0, & \quad \text{sell stock at } k. \end{aligned}$$

The portfolio values for $\{I(k;\omega)\}_{k=0}^T$ with $x = (x_0, x_1)$ are given by,

$$X^0(k;x,I) = x^0 - \sum_{\ell=0}^{k} h(I(\ell))\hat{P}(\ell)$$
$$X^1(k;x,I) = x^1 + \sum_{\ell=0}^{k} I(\ell).$$

Here

$$h(z) = \begin{cases} (1+\lambda_0)z, & z > 0 \\ (1-\lambda_1)z, & z \leq 0, \end{cases}$$

where $\lambda_1, \lambda_0 > 0$ are the proportional costs for selling and buying the stock, respectively.

We are interested in the following finance problems.

**Pricing derivative:** Let $Y = (Y^0, Y^1)$ be $\mathcal{F}_T$ measurable. We define $\pi^*(Y)$ the minimum of $x_0 P^0(0)$ such that for some $I$,

$$Y^0 \leq X^0(T;(x^0,0),I), \ Y^1 \leq X^1(T;(x^0,0),I).$$

We say $\pi^*(Y)$ is the price of $Y = (Y^0, Y^1)$. The problem is to calculate $\pi^*(Y)$. Another important problem is to obtain a strategy $I(\cdot)$ such that for $x^0 = \pi^*(Y)$,

$$Y^0 \leq X^0(T;(x^0,0),I), \ Y^1 \leq X^1(T;(x^0,0),I).$$

For the later use, we also define $\pi^*(Y;x_1)$ the minimum of $x_0 P^0(0)$ such that for some $I$,
$$Y^0 \leq X^0(T;(x^0,x_1),I), \ Y^1 \leq X^1(T;(x^0,x_1),I).$$
Then $\pi^*(Y;x_1) = \pi^*(\tilde{Y})$, where $\tilde{Y}^0 = Y^0, \tilde{Y}^1 = Y^1 - x_1$.

**Optimizing expected utility:** Let $U$ be a utility function. Let $(x^0, x^1)$ be given such that
$$x^0 P^0(0) - h(-x^1)P^1(0) > 0.$$
$V(x^0, x^1)$ is the maximum of
$$E[U(X^0(T; (x^0, x^1), I)P^0(T) - h(-X^1(T; (x^0, x^1), I))P^1(T))],$$
where $I(\cdot)$ is an admissible strategy: $k = 0, 1, 2, \ldots, T$,
$$X^0(k; (x^0, x^1), I)P^0(k) - h(-X^1(k; (x^0, x^1), I))P^1(k) \geq 0.$$
We want to calculate $V(x^0, x^1)$ and find a strategy $I$ that attains the maximum.

## 3. Price systems and a price formula

**Definition.** We say that $(\rho^0, \rho^1)$ is a price system if $\rho^0, \rho^1$ are positive random variables such that

(a) $E[\rho^0] = P^0(0)$;
(b) Define
$$R(k; \omega) = \frac{\rho^1(k; \omega)}{\rho^0(k; \omega)} \frac{1}{\hat{P}(k; \omega)},$$
$$\rho^0(k; \omega) = E[\rho^0 | \mathcal{F}_k], \quad \rho^1(k; \omega) = E[\rho^1 | \mathcal{F}_k].$$
Then
$$(1 - \lambda_1) \leq R(k; \omega) \leq (1 + \lambda_0), \quad k = 0, 1, 2, \ldots, T.$$
We denote $\mathcal{P}(\lambda_0, \lambda_1)$ the family of price systems.

**Remark.** Assume there is an equivalent martingale measure $Q$. Then $\mathcal{P}(\lambda_0, \lambda_1) \neq \emptyset$. In fact, define
$$\rho^0 = \frac{dQ}{dP} P^0(0)$$
$$\rho^1 = \rho^0 \hat{P}(T).$$
Then
$$\rho^0(k; \omega) = \frac{dQ}{dP}|_{\mathcal{F}_k} P^0(0)$$
$$\rho^1(k; \omega) = \rho^0(k; \omega)\hat{P}(k; \omega).$$
We can show that $(\rho^0, \rho^1)$ is a price system.

On the other hand, in the case $\lambda_0 = \lambda_1 = 0$, $(\rho^0, \rho^1)$ is a price system if and only if
$$\frac{dQ}{dP} = \rho^0(T)/P^0(0)$$
defines an equivalent martingale measure.

**Theorem 1** ([12]). *Assume $\mathcal{P}(\lambda_0, \lambda_1) \neq \emptyset$. Then*

(3.1) $$\pi^*(Y) = \sup_{\mathcal{P}(\lambda_0, \lambda_1)} E[Y^0 \rho^0 + Y^1 \rho^1].$$

**Remark.** If $\lambda_0 = \lambda_1 = 0$, then the above is

$$\pi^*(Y) = \sup_\rho E[\rho H \frac{P^0(0)}{P^0(T)}]$$

$$H = Y^0 P^0(T) + Y^1 P^1(T).$$

Here

$$\frac{dQ}{dP} = \rho$$

defines an equivalent martingale measure.

A similar result for diffusion models is given in [3].

## 4. Price system and optimal expected utility

In the following, we assume $P^0(k) = 1$ for all $k$.

Let $U$ be a strictly increasing utility function. Define

$$U^*(y) = \sup\{U(x) - xy; x \geq 0\}.$$

Define

$$V^*(\xi, x_1) = \inf\{E[U^*(\xi \rho^0(T))] + x_1 \xi E[\rho^1(T)]\}$$

**Theorem 2.** *We have*

$$V(x_0, x_1) \leq \inf_{\xi > 0}\{V^*(\xi, x_1) + x_0 \xi\}$$

*This is the same as*

(4.1) $$V(x_0, x_1) \leq \inf_{\xi > 0, \rho^0, \rho^1}\{E[U^*(\xi \rho^0(T))] + x_1 \xi E[\rho^1(T)] + x_0 \xi\}.$$

*Assume there is $\hat{\xi}, \hat{\rho}^0, \hat{\rho}^1$ that attains the infimum. Then the above equality holds. Moreover, there is an optimal strategy $\hat{I}$ for the portfolio optimization problem satisfying the following properties.*

(a) $X^0(T;(x_0,x_1),\hat{I}) = -U^{*'}(\hat{\xi}\hat{\rho}^0(T))$,
    $X^1(T;(x_0,x_1),\hat{I}) = 0$.
(b) $\hat{R}(l) = 1 + \lambda_0$ if $\hat{I}(l) > 0$,
    $\hat{R}(l) = 1 - \lambda_1$ if $\hat{I}(l) < 0$.

*Here $U^{*'}(\xi)$ denotes the derivative of $U^*(\xi)$.*

*Proof.* Let $I$ be a strategy.

(4.2) $$\begin{aligned}
U(X^0&(T;x,I) - h(-X^1(T;x,I))P^1(T)) \\
&\leq U^*(\xi \rho^0(T)) + \xi \rho^0(T)(X^0(T;x,I) - h(-X^1(T;x,I))P^1(T)) \\
&= U^*(\xi \rho^0(T)) + \xi(X^0(T;x,I)\rho^0(T) - h(-X^1(T;x,I))P^1(T)\rho^0(T)) \\
&\leq U^*(\xi \rho^0(T)) + \xi(X^0(T;x,I)\rho^0(T) + R(T)X^1(T;x,I)P^1(T)\rho^0(T)) \\
&= U^*(\xi \rho^0(T)) + \xi(X^0(T;x,I)\rho^0(T) + \rho^1(T)X^1(T;x,I)).
\end{aligned}$$

(4.3) $$\begin{aligned}
X^0(T;x,I)&\rho^0(T) + X^1(T;x,I)\rho^1(T) \\
&= x_0\rho^0(T) + x_1\rho^1(T) + (-\sum_{l=0}^{T} h(I(l))P^1(l)\rho^0(T) \\
&\quad + \sum_{l=0}^{T} I(l)\rho^1(T)).
\end{aligned}$$

$$
\begin{aligned}
& E[(-\sum_{l=0}^{T} h(I(l))P^1(l)\rho^0(T) + \sum_{l=0}^{T} I(l)\rho^1(T))] \\
& = \sum_{l=0}^{T} E[-h(I(l))P^1(l)\rho^0(l) + I(l)\rho^1(l)] \\
& = \sum_{l=0}^{T} E[(-h(I(l)) + R(l)I(l))P^1(l)\rho^0(l)] \\
& \leq 0.
\end{aligned}
\tag{4.4}
$$

Then we can deduce

$$
\begin{aligned}
& E[U(X^0(T;x,I) - h(-X^1(T;x,I))P^1(T))] \\
& \leq E[U^*(\xi\rho^0(T))] + \xi x_1 E[\rho^1(T)] + \xi x_0.
\end{aligned}
\tag{4.5}
$$

This is true for all $\rho^0, \rho^1$. The first result follows.

Assume $\hat{\xi}, \hat{\rho}^0, \hat{\rho}^1$ attains infimum in (4.1). Then

$$
E[U^{*'}(\hat{\xi}\hat{\rho}^0(T))\hat{\rho}^0(T)] + x_1 E[\hat{\rho}^1(T)] + x_0 = 0.
\tag{4.6}
$$

On the other hand, take any $(\rho^0, \rho^1)$ and $0 < \alpha < 1$, we have

$$
E[U^*(\hat{\xi}(\alpha\rho^0(T) + (1-\alpha)\hat{\rho}^0(T)))] + x_1\hat{\xi}E[\alpha\rho^1(T) + (1-\alpha)\hat{\rho}^1(T)] + x_0\hat{\xi}
$$

takes minimum at $\alpha = 0$. We have

$$
\begin{aligned}
& E[U^{*'}(\hat{\xi}\hat{\rho}^0(T))(\rho^0(T) - \hat{\rho}^0(T))] \\
& + x_1\hat{\xi}E[\rho^1(T) - \hat{\rho}^1(T)] \geq 0.
\end{aligned}
\tag{4.7}
$$

Take

$$
\hat{Y}^0 = -U^{*'}(\hat{\xi}\hat{\rho}^0(T)), \ \hat{Y}^1 = 0.
$$

(4.7) implies

$$
\pi^*(\hat{Y}; x_1) = E[-U^{*'}(\hat{\xi}\hat{\rho}^0(T))\hat{\rho}^0(T)] - x_1\hat{\xi}E[\hat{\rho}^1(T)] = x_0.
$$

Here we use (4.6) and Theorem 1.

By the definition of $\pi^*(\hat{Y}; x_1)$, there is a strategy $\hat{I}$ such that

$$
x_0 - \sum_{l=0}^{T} h(\hat{I}(l))P^1(l) \geq \hat{Y}^0,
$$

$$
x_1 + \sum_{l=0}^{T} \hat{I}(l) \geq \hat{Y}^1.
$$

Therefore,

$$
\begin{aligned}
& X^0(T; (x_0, x_1), \hat{I}) \geq -U^{*'}(\hat{\xi}\hat{\rho}^0(T)), \\
& X^1(T; (x_0, x_1), \hat{I}) \geq 0.
\end{aligned}
\tag{4.8}
$$

$$
\begin{aligned}
& U(X^0(T;x,\hat{I}) - h(-X^1(T;x,\hat{I}))P^1(T)) \\
& \geq U(-U^{*'}(\hat{\xi}\hat{\rho}^0(T))) \\
& = U^*(\hat{\xi}\hat{\rho}^0(T)) - \hat{\xi}\hat{\rho}^0(T)U^{*'}(\hat{\xi}\hat{\rho}^0(T)).
\end{aligned}
\tag{4.9}
$$

Then

$$(4.10) \quad \begin{aligned} & E[U(X^0(T;x,\hat{I})-h(-X^1(T;x,\hat{I}))P^1(T))] \\ & \geq E[U^*(\hat{\xi}\hat{\rho}^0(T))] - \hat{\xi}E[\hat{\rho}^0(T)U^{*\prime}(\hat{\xi}\hat{\rho}^0(T)))] \\ & = E[U^*(\hat{\xi}\hat{\rho}^0(T))] + \hat{\xi}(x_1 E[\hat{\rho}^1(T)] + x_0) \\ & \geq V(x_0, x_1). \end{aligned}$$

Therefore, by the definition of $V(x_0, x_1)$, the inequalities become equalities in the above relation. We see (a) follows from the equalities in (4.8),(4.9) and (4.10). On the other hand, (4.2),(4.3), (4.4) and (4.5) also become equalities for $I = \hat{I}$, then (b) follows. This completes the proof. □

## 5. Binomial model and price systems

We take $P_k^0 = 1$ for all $k$. $0 < d < 1 < u$, $\lambda_0, \lambda_1 > 0$.
  The sample space is given by

$$\Omega = \{(a_1, a_2, \ldots, a_T); a_i \in \{u, d\}\}.$$

For $\omega = (a_1, a_2, \ldots, a_T) \in \Omega$, denote

$$\omega^k = (a_1, a_2, \ldots, a_k).$$

The price of stock is

$$P_k^1(\omega) = P_0^1 a_1 a_2 \cdots a_k.$$

We also write $P_k^1(\omega) = P_k^1(\omega^k)$.
  $\mathcal{F}_k$ is the $\sigma$-algebra generated by $P_t^1, t \leq k$. A function defined on $\Omega$ measurable w.r.t. $\mathcal{F}_k$ is given by $f(\omega^k)$.
  For $\omega = (a_1, a_2, \ldots, a_T) \in \Omega$, the probability is given by

$$P(\{\omega\}) = p^m(1-p)^{T-m},$$

where $m$ is the number of $k$ such that $a_k = u$, $0 < p < 1$.
  $\rho^0(k), \rho^1(k)$ are given by

$$\rho^0(k) = E[\rho^0|\mathcal{F}_k], \quad \rho^1(k) = E[\rho^1|\mathcal{F}_k].$$

We have the characterization of $\rho^0(k), \rho^1(k)$:

$$(PS1) \quad \begin{aligned} \rho^0(k, \omega^k) &= p\rho^0(k+1, (\omega^k, u)) + (1-p)\rho^0(k+1, (\omega^k, d)), \\ \rho^1(k, \omega^k) &= p\rho^1(k+1, (\omega^k, u)) + (1-p)\rho^1(k+1, (\omega^k, d)). \end{aligned}$$

$$(PS2) \quad (1-\lambda_1)P_k^1 \leq \frac{\rho^1(k)}{\rho^0(k)} \leq (1+\lambda_0)P_k^1, k = 1, 2, 3, \ldots, T.$$

It is convenient to consider

$$A(k) = \frac{\rho^1(k)}{\rho^0(k)}, \quad k = 0, 1, \ldots, T.$$

We can now describe the price systems in a binomial market. We omit the easy proof.

**Theorem 3** (Binomial model). *Let $\omega = (a_1, a_2, \ldots, a_T) \in \Omega$. Given $A_0$ a positive constant such that*
$$(1 - \lambda_1)P_0^1 \leq A_0 \leq (1 + \lambda_0)P_0^1.$$

Denote $\rho^0(0) = 1, \rho^1(0) = \rho^0(0)A_0$. Take positive constants $A^u, A^d$ such that
$$\min\{A^u, A^d\} < A_0 < \max\{A^u, A^d\}$$

and
$$(1 - \lambda_1)P_0^1 u \leq A^u \leq (1 + \lambda_0)P_0^1 u,$$
$$(1 - \lambda_1)P_0^1 d \leq A^d \leq (1 + \lambda_0)P_0^1 d.$$

If $a_1 = u$
$$\rho^0(1) = \frac{1}{p}\frac{A_0 - A^d}{A^u - A^d},$$
$$A_1 = A^u.$$

If $a_1 = d$
$$\rho^0(1) = \frac{1}{1-p}\frac{A^u - A_0}{A^u - A^d},$$
$$A_1 = A^d.$$

Define
$$\rho^1(1) = \rho^0(1)A_1.$$

Assume we have defined $A_0, A_1, \ldots, A_k$ and
$$\rho^0(1), \rho^0(2), \ldots, \rho^0(k), \rho^1(1), \rho^1(2), \ldots, \rho^1(k).$$

Take $A^u, A^d$ measurable w.r.t. $\mathcal{F}_k$ such that
$$\min\{A^u, A^d\} < A_k < \max\{A^u, A^d\}$$

and
$$(1 - \lambda_1)P_k^1 u \leq A^u \leq (1 + \lambda_0)P_k^1 u,$$
$$(1 - \lambda_1)P_k^1 d \leq A^d \leq (1 + \lambda_0)P_k^1 d.$$

If $a_{k+1} = u$,
$$\rho^0(k+1) = \rho^0(k)\frac{1}{p}\frac{A_k - A^d}{A^u - A^d},$$
$$\rho^1(k+1) = \rho^0(k+1)A^u;$$

if $a_{k+1} = d$,
$$\rho^0(k+1) = \frac{1}{1-p}\frac{A^u - A_k}{A^u - A^d},$$
$$\rho^1(k+1) = \rho^0(k+1)A^d.$$

Then $\rho^0(k), \rho^1(k)$ satisfy (PS1) and (PS2).

## 6. Binomial model: control problems for pricing derivatives

Assume $Y = (Y^0, Y^1)$ is given by

$$Y^0 = Y^0(P_T^1), \ Y^1 = Y^1(P_T^1).$$

Then the price $\pi^*(Y)$ is given by

$$\pi^*(Y) = \sup_{\rho^0(T), \rho^1(T)} E[\rho^0(T)Y^0(P_T^1) + \rho^1(T)Y^1(P_T^1)].$$

This can be rewritten as

$$\pi^*(Y) = \sup_{\rho^0(T), \rho^1(T)} E[\rho^0(T)(Y^0(P_T^1) + A_T Y^1(P_T^1))].$$

with $A_k$ and $\rho^0(k)$ described in Theorem 3. This is viewed as a stochastic control problem. The state variables are given by $P_k^1, A_k, \rho^0(k)$ and the control variables are $A_k^u, A_k^d$.

The dynamical programming can be described as follows.

For $S > 0$ and $A$ satisfying

$$(1 - \lambda_1)S \leq A \leq (1 + \lambda_0)S,$$

define

$$W_k(S, A) = \sup E[\frac{\rho^0(T)}{\rho^0(k)}(Y^0(P_T^1) + A_T Y^1(P_T^1)) | P_k^1 = S, A_k = A]$$

Then

$$\pi^*(Y) = \sup_{(1-\lambda_1)S \leq A \leq (1+\lambda_0)S} W_0(S, A)$$

for $P_0^1 = S$. And for $0 \leq k < l \leq T$,

$$W_k(S, A) = \sup E[\frac{\rho^0(l)}{\rho^0(k)} W_l(P_l^1, A_l) | P_k^1 = S, A_k = A]$$

It follows a recursive scheme backward in time.

(D1) $W_T(S, A) = Y^0(S) + AY^1(S)$;
(D2) For $(1 - \lambda_1)S \leq A \leq (1 + \lambda_0)S$,

$$W_k(S, A) = \sup\{\frac{A - A^d}{A^u - A^d} W_{k+1}(Su, A^u) + \frac{A^u - A}{A^u - A^d} W_{k+1}(Sd, A^d)\},$$

the maximization is taken over

$$\min\{A^u, A^d\} < A < \max\{A^u, A^d\}$$
$$(1 - \lambda_1)Su \leq A^u \leq (1 + \lambda_0)Su,$$
$$(1 - \lambda_1)Sd \leq A^d \leq (1 + \lambda_0)Sd.$$

(D3) For $P_0^1 = S$,

$$\pi^*(Y) = \sup_{(1-\lambda_1)S \leq A \leq (1+\lambda_0)S} \{W_0(S, A)\}.$$

It can be restated as follows.

**Theorem 4.** *We have*

$$W_k(S, A) = \sup\{\alpha W_{k+1}(Su, A^u) + (1-\alpha)W_{k+1}(Sd, A^d)\},$$

*the maximization is taken over*

$$0 < \alpha < 1,$$
$$\alpha A^u + (1-\alpha)A^d = A,$$

*and*

$$(1-\lambda_1)Su \leq A^u \leq (1+\lambda_0)Su,$$
$$(1-\lambda_1)Sd \leq A^d \leq (1+\lambda_0)Sd.$$

$W_k(S,A)$ *is piecewise linear in* $A$ *for all* $S > 0$ *and* $k = 0, 1, \ldots, T$.

The main questions consist of the following. How to calculate $W_k(S, A)$? How to obtain an optimal strategy to super hedge $Y$ from $W_k(S,A)$? Some answers can be found in [6].

## 7. Binomial model: optimizing expected utility and control problem

We take

$$U(x) = \frac{1}{\gamma}x^\gamma, 0 < \gamma < 1.$$

Then

$$U^*(\xi) = -\frac{1}{\mu}\xi^\mu,$$
$$\mu = \frac{\gamma}{\gamma - 1}.$$

Then

(7.1) $$V(x_0, x_1) = \inf_{\xi > 0}\{V^*(\xi, x_1) + x_0\xi\},$$

$$V^*(\xi, x_1) = \inf\{-\frac{1}{\mu}\xi^\mu E[(\rho^0(T))^\mu] + \xi x_1 E[\rho^1(T)]\}.$$

We shall consider

$$-\frac{1}{\mu}\xi^\mu E[(\rho^0(T))^\mu] + \xi x_1 E[\rho^1(T)]$$

conditioning on $P^1(0) = S, A(0) = A$. This is equal to

$$-\frac{1}{\mu}\xi^\mu E[(\rho^0(T))^\mu] + \xi x_1 A.$$

We consider

$$V_k(S, A) = \inf E[(\frac{\rho^0(T)}{\rho^0(k)})^\mu | \mathcal{F}_k].$$

The follwoing is an iterative scheme to calculate $V_k(S, A), k = 0, 1, \ldots$.
(PD1) $V_T(S, A) = 1$,
(PD2) $k = 0, 1, 2, \ldots$,

(7.2) $$V_k(S, A) = \inf\{p^{1-\mu}(\frac{A - A^d}{A^u - A^d})^\mu V_{k+1}(Su, A^u) \\ + (1-p)^{1-\mu}(\frac{A^u - A}{A^u - A^d})^\mu V_{k+1}(Sd, A^d)\}$$

The inf is taken over the $A^d, A^u$ satisfying

$$\min\{A^u, A^d\} < A < \max\{A^u, A^d\}$$

(7.3)
$$(1 - \lambda_1)Su \le A^u \le (1 + \lambda_0)Su,$$
$$(1 - \lambda_1)Sd \le A^d \le (1 + \lambda)Sd.$$

(7.2) can be reformulated as follows.

(7.2)' $\quad V_k(S, A) = \inf\{p^{1-\mu}\alpha^\mu V_{k+1}(Su, A^u) + (1-p)^{1-\mu}(1-\alpha)^\mu V_{k+1}(Sd, A^d)\}$

where $0 \le \alpha \le 1$ and (7.3) and (7.4) hold,

(7.4) $$\alpha A^u + (1-\alpha)A^d = A.$$

We consider $V_{T-1}(S, A)$:

$$V_{T-1}(S, A) = \inf\{p^{1-\mu}\alpha^\mu + (1-p)^{1-\mu}(1-\alpha)^\mu\},$$

where (7.3), (7.4) hold. Denote $\hat{V}_{T-1}(A) = V_{T-1}(S, SA)$. For

$$(1 - \lambda_1) \le A \le (1 + \lambda_0),$$

$$\hat{V}_{T-1}(A) = \inf\{p^{1-\mu}\alpha^\mu + (1-p)^{1-\mu}(1-\alpha)^\mu\},$$

there are $A^u, A^d$ such that

(7.3)'
$$(1 - \lambda_1) \le A^u \le (1 + \lambda_0),$$
$$(1 - \lambda_1) \le A^d \le (1 + \lambda_0).$$

(7.4)' $$\alpha u A^u + (1-\alpha)dA^d = A,$$

In general,

$$\hat{V}_k(A) = V_k(S, SA), \ (1 - \lambda_1) \le A \le (1 + \lambda_0).$$

Then

$$\hat{V}_k(A) = \inf\{p^{1-\mu}\alpha^\mu \hat{V}_{k+1}(A^u) + (1-p)^{1-\mu}(1-\alpha)^\mu \hat{V}_{k+1}(A^d)\},$$

$0 \le \alpha \le 1$ satisfies (7.3)', (7.4)'. From

$$\hat{V}_k(A), \ (1 - \lambda_1) \le A \le (1 + \lambda_0),$$

we have

$$V_k(S, A) = \hat{V}_k(\frac{A}{S}), \ (1 - \lambda_1)S \le A \le (1 + \lambda_0)S.$$

**Theorem 5.** *Assume $x_0 - h(-x_1)S > 0$. Then for $P^1(0) = S$,*

$$V(x_0, x_1) = \frac{1}{\gamma}\inf\{(x_0 + x_1SR)^\gamma (\hat{V}_0(R))^{1-\gamma}\}$$

*the infimum is taken over $(1 - \lambda_1) \le R \le (1 + \lambda_0)$.*
*In particular, if $1 < pu + (1-p)d$ and $x_1 \le 0$, then*

$$V(x_0, x_1) = \frac{1}{\gamma}(x_0 + x_1S(1 + \lambda_0))^\gamma (\hat{V}_0(1 + \lambda_0))^{1-\gamma}.$$

*If $1 > pu + (1-p)d$ and $x_1 \ge 0$, then*

$$V(x_0, x_1) = \frac{1}{\gamma}(x_0 + x_1S(1 - \lambda_1))^\gamma (\hat{V}_0(1 - \lambda_1))^{1-\gamma}.$$

*Proof.* By (7.1)

$$V(x_0, x_1) = \inf\{-\frac{1}{\mu}\xi^\mu \hat{V}_0(R) + \xi x_1 SR + \xi x_0\},$$

the *inf* is taken over $\xi > 0, (1 - \lambda_1) \leq R \leq (1 + \lambda_0)$.

$$\hat{\xi} = (\frac{1}{\hat{V}_0(R)}(x_0 + x_1 SR))^{\frac{1}{\mu-1}}$$

takes minimum. The rest follows from this and Theorem 6 below. □

**Theorem 6.** *Assume* $1 < pu + (1-p)d$. *Then for*

$$(1 - \lambda_1)(pu + (1-p)d)^{T-k} \leq A \leq (1 + \lambda_0),$$

$\hat{V}_k(A) = 1$. *For other* $A$, $\hat{V}_k(A) > 1$ *and is decreasing in* $A$.
*Assume* $pu + (1-p)d < 1$. *Then for*

$$(1 - \lambda_1)S \leq A \leq (1 + \lambda_0)S(pu + (1-p)d)^{T-k},$$

$\hat{V}_k(A) = 1$. *For other* $A$, $\hat{V}_k(A) > 1$ *and is increasing in* $A$.
$\hat{V}_k(A)$ *is nonincreasing in* $k$ *for fixed* $A$.

*Proof.* We only consider $1 < pu + (1-p)d$. Define

$$f(\alpha) = p^{1-\mu}\alpha^\mu + (1-p)^{1-\mu}(1-\alpha)^\mu, \ 0 < \alpha < 1.$$

$f$ takes minimum at $\alpha = p$, $f(p) = 1$ and $f$ is decreasing on $(0, p]$ and increasing on $[p, 1)$.
Given $A$,

$$(1 - \lambda_1) \leq A \leq (1 + \lambda_0).$$

We consider

$$\inf\{f(\alpha)\}.$$

The infimum is taken over $\alpha$ such that there are $A^u, A^d$ satisfying

$$\alpha u A^u + (1-\alpha)d A^d = A,$$

and

$$(1 - \lambda_1) \leq A^u \leq (1 + \lambda_0),$$
$$(1 - \lambda_1) \leq A^d \leq (1 + \lambda_0).$$

We consider the cases,

(i) $(1 - \lambda_1) \leq A \leq (1 - \lambda_1)u$;
(ii) $(1 - \lambda_1)u \leq A \leq (1 + \lambda_0)d$;
(iii) $(1 + \lambda_0)d \leq A \leq (1 + \lambda_0)$.

Assume (i),

$$\alpha = \frac{A - dA^d}{uA^u - dA^d}.$$

For each $(1 - \lambda_1) \leq A^u \leq (1 + \lambda_0)$, the range of $\alpha$ defined above taken over

$$(1 - \lambda_1)d \leq dA^d \leq A$$

is $[0, (A - (1 - \lambda_1)d)/(uA^u - (1 - \lambda_1)d)$. Take the union of these sets over all

$$(1 - \lambda_1) \leq A^u \leq (1 + \lambda_0),$$

we have $[0, (A - (1 - \lambda_1)d)/(1 - \lambda_1)(u - d)]$. If $p$ is in this interval, then $\hat{V}_{T-1}(A) = 1$. The condition $p$ is in this interval is the same as

$$A \geq (1 - \lambda_1)(pu + (1 - p)d).$$

Therefore,

$$\hat{V}_{T-1}(A) = 1, (1 - \lambda_1)(pu + (1 - p)d) \leq A \leq (1 - \lambda_1)u.$$

On the other hand, if

$$(1 - \lambda_1) \leq A \leq (1 - \lambda_1)(pu + (1 - p)d),$$

the infimum of $f(\alpha)$ on

$$[0, (A - (1 - \lambda_1)d)/(1 - \lambda_1)(u - d)]$$

is

$$f(\frac{A - (1 - \lambda_1)d}{(1 - \lambda_1)(u - d)}).$$

Therefore,

$$\hat{V}_{T-1}(A) = f(\frac{A - (1 - \lambda_1)d}{(1 - \lambda_1)(u - d)})$$

if

$$(1 - \lambda_1) \leq A \leq (1 - \lambda_1)(pu + (1 - p)d).$$

Assume (ii). We consider $A \leq uA^u \leq (1 + \lambda_0)u$. The range of $\alpha$ is given by $[0, 1]$. Therefore, $\hat{V}_{T-1}(A) = 1$.

Assume (iii). For each $A \leq uA^u \leq (1 + \lambda_0)u$, the range of $\alpha$ of

$$(1 - \lambda_1) \leq A^d \leq (1 + \lambda_0)$$

is

$$[\frac{A - (1 + \lambda_0)d}{uA^u - (1 + \lambda_0)d}, \frac{A - (1 - \lambda_1)d}{uA^u - (1 - \lambda_1)d}].$$

Take the union of these sets over all $A^u$ gives

$$[\frac{A - (1 + \lambda_0)d}{(1 + \lambda_0)(u - d)}, 1].$$

We can check $p$ is in this set. Then $\hat{V}_{T-1}(A) = 1$.
We conclude

$$\hat{V}_{T-1}(A) = f(\frac{A - (1 - \lambda_1)d}{(1 - \lambda_1)(u - d)})$$

if

$$(1 - \lambda_1) \leq A \leq (1 - \lambda_1)(pu + (1 - p)d),$$

and

$$\hat{V}_{T-1}(A) = 1$$

if
$$(1 - \lambda_1)(pu + (1-p)d) \leq A \leq (1+\lambda_0).$$

$\hat{V}_{T-1}(A)$ is decreasing in $A$.

We can continue this argument for other $\hat{V}_k(A)$ to prove that $\hat{V}_k(A) = 1$ if
$$(1 - \lambda_1)(pu + (1-p)d)^{T-k} \leq A \leq (1+\lambda_0),$$

and for other $A$, $\hat{V}_k(A) > 1$. To prove the nonincreasing of $\hat{V}_k(A)$ in $A$ needs additional argument. We have the following observation. Let $g$ be nonincreasing. Consider
$$\hat{g}(A) = \inf\{p^{1-\mu}\alpha^\mu g(A^u) + (1-p)^{1-\mu}(1-\alpha)^\mu g(A^d)\},$$

where the "inf" is taken over $0 < \alpha < 1$ and $A^u, A^d$ satisfying $(7.3)'$ and $(7.4)'$. We define
$$\bar{g}(A) = g(A), \ (1 - \lambda_1) \leq A \leq (1+\lambda_0),$$
$$\bar{g}(A) = \infty, \ A < (1 - \lambda_1),$$
$$\bar{g}(A) = g((1+\lambda_0)), \ A > (1+\lambda_0).$$

We claim

(7.5) $$\hat{g}(A) = \inf\{p^{1-\mu}\alpha^\mu \bar{g}(A^u) + (1-p)^{1-\mu}(1-\alpha)^\mu \bar{g}(A^d)\},$$

where the "inf" is taken over $0 < \alpha < 1$ and $A^u, A^d$ satisfying $(7.4)'$. First, it is easy to see that the quantity defined by the righthand side of (7.5) is not smaller than $\hat{g}(A)$. To prove the opposite inequality, we observe that for a given $0 < \alpha < 1$ and $A^u, A^d$ satisfying $(7.4)'$, if $(7.3)'$ does not hold, says
$$A^u > (1+\lambda_0).$$

We define $\bar{A}^u = (1+\lambda_0)$ and $\bar{A}^d$ by the relation,
$$\alpha u(1+\lambda_0) + (1-\alpha)d\bar{A}^d = A.$$

Then $\bar{A}^d > A^d$. We see $\bar{A}^u, \bar{A}^d$ satisty $(7.3)'$ and $(7.4)'$ and
$$p^{1-\mu}\alpha^\mu \bar{g}(A^u) + (1-p)^{1-\mu}(1-\alpha)^\mu \bar{g}(A^d)$$
$$\geq p^{1-\mu}\alpha^\mu g(\bar{A}^u) + (1-p)^{1-\mu}(1-\alpha)^\mu g(\bar{A}^d)$$

by the property that $g$ is nonincreasing. Using this observation, we can deduce that the quantity defined by the righthand side of (7.5) is not smaller than $\hat{g}(A)$.

Now from (7.5) it is easy to see that $\hat{g}$ is nonincreasing. In fact, let $B = \lambda A > 0$ for a $\lambda > 1$. Let $0 < \alpha < 1$ and $A^u, A^d > 0$ satisfying
$$\alpha u A^u + (1-\alpha)dA^d = A.$$

We take $B^u = A^u\lambda, B^d = A^d\lambda$. Then
$$\alpha u B^u + (1-\alpha)dB^d = B.$$

We have
$$p^{1-\mu}\alpha^\mu \bar{g}(A^u) + (1-p)^{1-\mu}(1-\alpha)^\mu \bar{g}(A^d)$$
$$\geq p^{1-\mu}\alpha^\mu g(B^u) + (1-p)^{1-\mu}(1-\alpha)^\mu g(B^d)$$
$$\geq \hat{g}(B).$$

This is true for any $\alpha, A^u, A^d$. Therefore, $\hat{g}(A) \geq \hat{g}(B)$.

Finally, we denote $\hat{g}(A) = Hg(A)$. Then $H$ has the property that $g_1(A) \geq g_2(A)$ for all $A$ implies $Hg_1(A) \geq Hg_2(A)$ for all $A$. Take $g = 1$. Then

$$Hg = \hat{V}_{T-1}.$$

We have proved $\hat{V}_{T-1} \geq 1$. That is,

$$Hg \geq g.$$

We note

$$\hat{V}_k = H\hat{V}_{k+1}.$$

From these, by induction, we can show $\hat{V}_k \geq \hat{V}_{k+1}$. This completes the proof. □

**Corollary 7.** *Assume $1 < pu + (1-p)d$ and*

$$(1 - \lambda_1)(pu + (1-p)d)^T \leq (1 + \lambda_0).$$

*If $x_1 \leq 0$, then buy-and-hold is an optimal strategy.*

*Similarly, assume $1 > pu + (1-p)d$ and*

$$(1 + \lambda_0)(pu + (1-p)d)^T \geq (1 - \lambda_1).$$

*If $x_1 \geq 0$, then sell-and-hold is an optimal strategy.*

*Proof.* Assume $1 < pu + (1-p)d$ and $x_1 \leq 0$. From Theorem 5 and 6,

$$V(x_0, x_1) = \frac{1}{\gamma}(x_0 + x_1 S(1 + \lambda_0))^\gamma.$$

Buy-and-hold achives this value and hence is an optimal strategy. Other result can be proved similarly. □